\documentclass[12pt,leqno]{article}
\usepackage{amsfonts, amssymb}
\usepackage{tikz}
\usetikzlibrary{snakes}

\newtheorem{theorem}{Theorem}
\newtheorem{lemma}{Lemma}

\begin{document}

\begin{center}

{\large On a Certain Function with Negative Coefficients }

\vspace{10pt}

Dov Aharonov

\vspace{10pt}

\end{center}
\vspace{20pt}
\begin{abstract}

\noindent In what follows we improve an inequality
related to matrix theory.\\
T. Laffey  proved (2013) a weaker form of this inequality \cite{2}.

\end{abstract}

\vspace{20pt}

\begin{theorem}
Given 
$\rho = \sum\limits_{k=1}^{n}\mu_k \leq 1$ ,
$\mu_k > 0$ ,
$k=1,2,\ldots,n$,

$c_1 >c_2>\cdots >c_n>0$

Then for $n \geq 2$ , or $n=1$ $\mu_1<1$:

\begin{equation}
\sum\limits_{j=1}^{\infty} \pi (1-c_jt)^{\mu_j} = 1-\sum\limits_{j=1}^{\infty}D_j t^j
\ \mbox{we have} \
D_j >0 \ , \ j=1,2,\ldots
\end{equation}
\end{theorem}

\begin{lemma}
In order to prove the theorem it is enough to prove it for $\rho =1$.
\end{lemma}

\noindent {\bf Proof of Lemma 1}

Assume that $\sum\limits_{j=1}^{n}\mu_j=1$ and that for this case the theorem is proved. Then for $0 < \rho <1$ we have

\[
(1-\sum D_jt^j)^{\rho}=
1-\rho\sum D_jt^j + \rho(\rho-1)\left(\sum D_jt^j\right)^{2}-
\frac{1}{6}\rho(\rho-1)(\rho-2)\left(\sum D_jt^j\right)^{3}+ \cdots
\]
\[=1-\sum\limits_{j=1}^{\infty}\tilde{D}_j t^j\]
Hence $\tilde{D}_j > 0$ 

\begin{flushright}
$\blacksquare$
\end{flushright}

From now on we may assume, without loss of generality,
\begin{equation}
\sum\limits_{j=1}^{n}\mu_j =1
\end{equation}
\begin{lemma}
Assuming $(1)$ and $(2)$ we have
\begin{equation}
\sum\limits_{j=1}^{n}c_j \mu_j= D_1 \ , \ D_1 > c_n
\end{equation}
\end{lemma}

\noindent {\bf Proof of Lemma 2}

\[
\prod\limits_{j=1}^{n}(1-c_jt)^{\mu_j} =
\prod\limits_{j=1}^{n}(1-\mu_jc_jt+O(t^2)) =
1-\left(\sum\limits_{j=1}^{n}\mu_jc_j\right)t+ O(t^2)=
1-D_1t+O(t^2)
\]
Hence, the first part is confirmed.

For the second part recall the monotonicity of $c_j$.

We have

\[
D_1 = \sum\limits_{j=1}^{n}c_j\mu_j > c_n 
\sum\limits_{j=1}^{n} \mu_j = c_n
\]

\begin{flushright}
$\blacksquare$
\end{flushright}

The proof of our theorem will be by induction. For $n=1,2$ the proof follows very simply. Details are omitted.

Now assume that the theorem is proved for $n-1$ ($n \geq 3$).
We want to prove it for $n$.
Suppose that this is not the case.

It will be convenient to use the notation
\begin{equation}
L_k = \prod\limits_{j=1}^{k}(1-c_jt)^{\mu_j}
\end{equation}

Obviously, $L_k = L_{k-1}(1-c_kt)^{\mu_k}$. In particular,
\[
L_n = L_{n-1}(1-c_nt)^{\mu_n}=1-\sum\limits_{j=1}^{\infty}D_jt^j
\]
By our assumption, there exists among the coefficients at least one which is not negative. Take the smallest index, say $m$, for which this is true. Hence
\begin{equation}
-D_j < 0 \ , \ 1 \leq j \leq m-1 \ , \ -D_m \geq 0
\end{equation}
Time has come to use a simple idea, but useful.

From $(4)$
\[
L_n = L_{n-2}(1-c_{n-1}t)^{\mu_{n-1}}(1-c_nt)^{\mu_n}
\]
Then for $c_n=0$ $L_n=L_{n-1}$ and we are back to the case of $n-1$ factors, i.e. the theorem is assumed to be correct by the induction assumption. The same follows if we assume $c_{n-1}=c_n$.

Indeed $(1-c_{n-1}t)^{\mu_{n-1}}(1-c_nt)^{\mu_n}=
(1-c_{n-1}t)^{\mu_{n-1}+\mu_n}$,
and again we have $n-1$ factors. (Note that if we want to prove a weaker result for $\mu_1=\mu_2=\cdots=\mu_n = \frac{1}{n}$ the proof does not work!!)

\begin{minipage}{3in}
\begin{center}

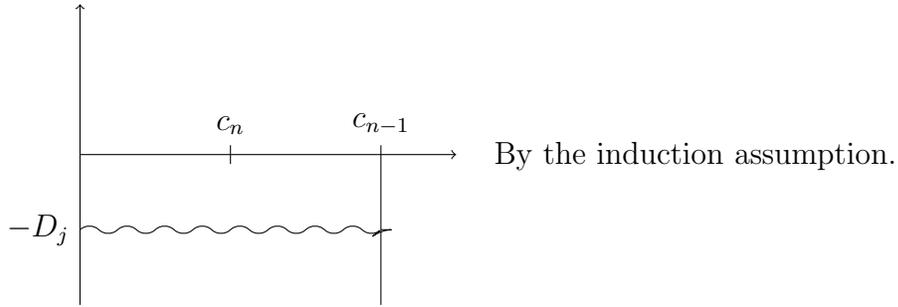
{\bf Figure 1}

\vspace{10pt}

\begin{tikzpicture}[scale=0.5]
\draw[->] (0,0)--(0,8);
\draw[->] (0,4)--(10,4);
\draw (8,0)--(8,4.25);
\draw (4,3.75)--(4,4.25);
\draw[snake,segment amplitude = 1mm, segment length = 5mm,rounded corners] (0,2)--(8,2);
\node[above] at (4,4.25) {$c_n$};
\node[above] at (8,4.25) {$c_{n-1}$};
\node[left] at (0,2) {$-D_j$};
\node at (0,-1) {$j=1,2,\ldots, n-1$};
\end{tikzpicture}
\end{center}
\end{minipage}\hspace{2pt}\begin{minipage}{2.5in}
By the induction assumption.
\end{minipage}

\vspace{10pt}

\begin{minipage}{3in}
\begin{center}

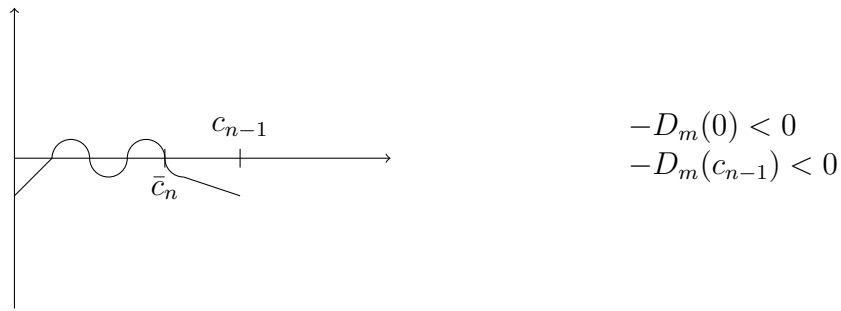
{\bf Figure 2}

\vspace{10pt}

\begin{tikzpicture}[scale=0.5]
\draw[->] (0,0)--(0,8);
\draw[->] (0,4)--(10,4);
\node[below] at (4,3.8) {$\bar{c}_n$};
\draw (4,3.75)--(4,4.25);
\node[above] at (6,4.25) {$c_{n-1}$};
\draw (6,3.75)--(6,4.25);
\draw (3,4) arc (180:0:.5);
\draw (2,4) arc (-180:0:.5);
\draw (1,4) arc (180:0:.5);
\draw (4,4) arc (-180:-90:.5);
\draw (4.5,3.5)--(6,3);
\draw (0,3)--(.5,3.5);
\draw (.5,3.5)--(1,4);
\end{tikzpicture}
\end{center}
\end{minipage}\hspace{2pt}\begin{minipage}{2.5in}
\[
\begin{array}{l}
-D_m(0) <0\\
-D_m(c_{n-1}) <0
\end{array}
\]
\end{minipage}

\vspace{10pt}

Consider now $c_n$ as a variable and $c_1, c_2, \ldots, c_{n-1}$ as fixed.
\[
0 \leq c_n \leq c_{n-1}
\]
Since $-D_m \geq 0$ at some point, say $\bar{c}_n$, for this interval and also
$-D_m(0) < 0$, $-D_m(c_{n-1}) < 0$ it follows by the mean value theorem that for some point in this interval, say 
$c^*_n$: $-D_m(c^*_n)=0$.

But there are only finite number of zeros of $-D_m$. This is true due to the fact that it is a polynomial function of the variable $c_n$. Thus, without loss of generality, $c^*_n$ may be taken as the largest zero in the interval $(0,c_{n-1})$.

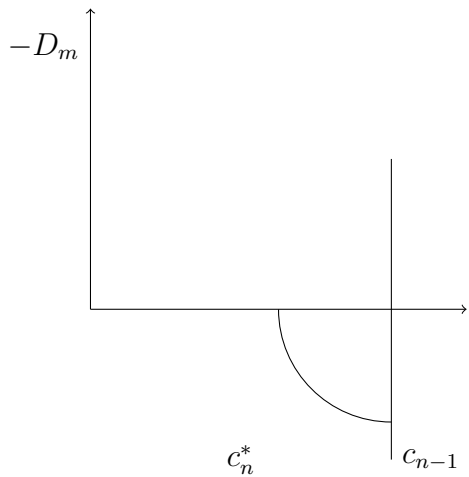
\begin{center}

{\bf Figure 3}

\vspace{10pt}

\begin{tikzpicture}[scale=0.5]
\draw[->] (0,0)--(0,8);
\draw[->] (0,4)--(10,4);
\draw (7,6)--(7,2);
\node[left] at (0,7) {$-D_m$}; 
\draw (0,2)-- (1,4);

\draw (7,3) arc (270:180:1);
\draw[fill] (6,4) circle (0.07);
\draw[snake,segment amplitude = 3mm, segment length = 10mm,rounded corners] (1,4)--(6,4);
\node at (6,5) {$c^*_n$}; 
\node[below] at (7,2) {$c_{n-1}$}; 
\end{tikzpicture}
\end{center}

\vspace{10pt}

\begin{center}

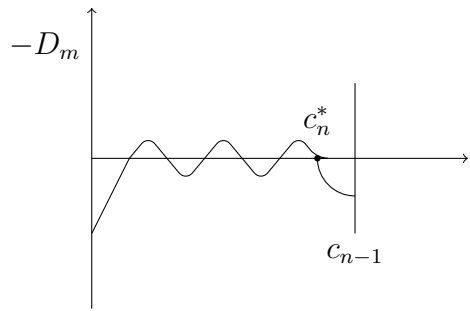
{\bf Figure 4}

\vspace{10pt}

\begin{tikzpicture}[scale=0.5]
\draw[->] (0,0)--(0,8);
\draw[->] (0,0)--(10,0);
\draw (8,-4)--(8,4);
\draw (8,-3) arc (270:180:3);
\node[left] at (0,7) {$-D_m$}; 
\node[right] at (8,-4) {$c_{n-1}$};
\node at (4,-4) {$c_n^*$};
\end{tikzpicture}
\end{center}

Thus (due to the fact that there are finite number of zeros as explained above) we have that the interval $(c^*_n,c_{n-1})$ is free of zeros.


\vspace{10pt}

In what follows we show a contradiction. this will be done by considering a sufficiently small interval $(c^*_n,c_{n}^*+ \varepsilon)$ for which $-D_m$ is
changed from zero at $c_n^*$ to a positive value. 

We now present the calculations leading to this assertion.

Indeed
\[
L_n(c_1,c_2,\ldots,c_{n-1},c^*_{n}+\varepsilon)=
L_{n-1}(c_1,c_2,\ldots,c_{n-1})(1-(c^*_{n}+\varepsilon)t)^{\mu_n}=
\]
\[
= L_{n-1}(1-c^*_nt-\varepsilon t)^{\mu_n}
\]
Also:

\[
(1-c_n^*t-\varepsilon t)^{\mu_n}=
(1-c_n^*t)^{\mu_n}\left(1-\frac{\varepsilon t}{1-c^*_nt}\right)^{\mu_n}=
(1-c_n^*t)^{\mu_n}\left(1-\frac{\varepsilon t}{1-c^*_nt} \mu_n +O(\varepsilon^2)\right)
\]
Thus
\[
L_n(c_1,c_2, \ldots ,c_{n-1},c^*_{n}+\varepsilon)=
[L_{n-1}(c_1,c_2, \ldots ,c_{n-1})(1-c_n^*t)^{\mu_n}]
\left(1-\frac{\varepsilon t}{1-c^*_nt} \mu_n +O(\varepsilon^2)\right)
\]
\[=(1-D_1 t - D_2t^2-\cdots -D_{m-1}t^{m-1}+0-D_{m+1}t^{m+1}+
\cdots)
\left(1-\frac{\varepsilon t\mu_n}{1-c^*_nt}+O(\varepsilon^2)\right)
=R \ \mbox{(notation)}
\]
To continue it will be convenient to use the notation
$\{F\}_K$ for the $k$-th coefficient of $F$.

Also denote $\mu_n\varepsilon$ by $\varepsilon^*$.
\[
\{R\}_m = \left\{\left(1-\sum\limits_{j=1}^{m-1}D_j t^j\right)\left(1-\frac{\varepsilon^* t}{1-c^*_nt}\right)\right\}_m +O(\varepsilon^2)
\]
Positive terms do not destroy positivity if ignored, hence
\[
\{R\}_m = \left\{-\frac{\varepsilon^*t}{1-c^*_nt}\right\}_m+
\left\{(D_1t)\left(\frac{\varepsilon^*t}{1-c^*_nt}\right)\right\}_m+ \ \mbox{positive terms} \ +O(\varepsilon^2)
\]
\[
\begin{array}{rl}
\{R\}_m \geq& -\varepsilon(c_n^*)^{m-1}+
(\varepsilon^*)\left\{D_1t^2\left(1+\sum\limits_{k=1}^{\infty}
(c_n^*t)^k\right)\right\}_m+O(\varepsilon^2)
\\
\{R\}_m \geq& -\varepsilon(c_n^*)^{m-1}+
\varepsilon^*D_1(c_n^*)^{m-2}+O(\varepsilon^*)
\\
\vspace{-10pt}\\
\{R\}_m \geq& \varepsilon^*(c_n^*)^{m-2}[-c_n^*+D_1]+O(\varepsilon^*)
\end{array}
\]
For $\varepsilon$ small enough this is positive by $(3)$ applied for
the special value $c_n^*$ of $c_n$.

Hence, positivity of the $n$-th coefficient is established provided $\varepsilon$ is small enough.
Thus we arrived at a contradiction, which ends the proof. $\blacksquare$

\vspace{10pt}

\begin{minipage}{3in}

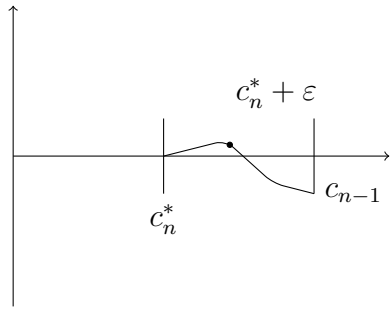
{\bf Figure 5}

\vspace{10pt}

\begin{center}
\begin{tikzpicture}[scale=0.5]
\draw[->] (0,0)--(0,8);
\draw[->] (0,4)--(10,4);
\draw (4,5)--(4,3);
\node[below] at (4,3) {$c^*_n$};
\draw[snake,segment amplitude = 4mm, segment length = 30mm,rounded corners] (4,4)--(8,3);
\draw (8,5)--(8,3);
\node[right] at (8,3) {$c_{n-1}$};
\draw[fill] (5.76,4.3) circle (.07);
\node[above] at (7,5) {$c^*_n+\varepsilon$};
\end{tikzpicture}
\end{center}
\end{minipage}\begin{minipage}{2in}
Contradiction

Compare with Figure 4

\end{minipage}

\vspace{1in}

\noindent\underline{\bf Acknowledgment}\\
I want to thank Prof Raphi Loewy for bringing to my attention the papers
\cite{1} and \cite{2}.

\vspace{0.5in}
\noindent
Department of mathematics-I.I.T,Haifa 32000.Israel.\\
E-mail address:\ \ \ dova@tx.technion.ac.il
\end{document}